\documentclass[leqno,12pt]{article}

\usepackage{amssymb, ulem}
\usepackage{euscript}
\usepackage{mathrsfs} 
\usepackage[latin1]{inputenc}
\usepackage{amsmath}
\usepackage{amsthm}
\usepackage{amsfonts}
\usepackage{enumerate}

\usepackage[colorlinks = true, linkcolor = blue, urlcolor  = blue, citecolor = blue, anchorcolor = blue]{hyperref} 

\usepackage{graphicx} 
\usepackage{float} 

\usepackage{ulem}

 \usepackage{pgf,tikz}
\usepackage{mathrsfs}
\usetikzlibrary{arrows}
\usetikzlibrary{shapes}
\usetikzlibrary{backgrounds}

\evensidemargin\oddsidemargin

\def\thesection{\arabic{section}}
\renewcommand{\theequation}{\thesection.\arabic{equation}}

\newtheorem{theorem}{Theorem}[section]

\newtheorem{proposition}[theorem]{Proposition}
\newtheorem{corollary}[theorem]{Corollary}

\newtheorem{notation}[theorem]{Notation}

\theoremstyle{definition}   
\newtheorem{remark}[theorem]{Remark}

\newcommand{\eqnsection}{
\renewcommand{\theequation}{\thesection.\arabic{equation}}
    \makeatletter
    \csname  @addtoreset\endcsname{equation}{section}
    \makeatother}
\eqnsection

\def\r{{\mathbb R}}
\def\e{{\mathbb E}}

\def\d{\mathrm{d}}

\begin{document}


\vglue20pt

\centerline{\large\bf  An infinite-dimensional representation of}
\centerline{\large\bf  the Ray--Knight theorems}

\bigskip
\bigskip

\centerline{by}

\medskip

\centerline{Elie A\"{i}d\'ekon\footnote{\scriptsize LPSM, Sorbonne Universit\'e Paris VI, NYU Shanghai, and Institut Universitaire de France, {\tt elie.aidekon@upmc.fr}},   Yueyun Hu\footnote{\scriptsize LAGA, Universit\'e Sorbonne Paris Nord, France, {\tt yueyun@math.univ-paris13.fr}}, 
and Zhan Shi\footnote{\scriptsize LPSM, Sorbonne Universit\'e Paris VI, France, {\tt zhan.shi@upmc.fr}}}

\bigskip
\bigskip

{\leftskip=2truecm \rightskip=2truecm \baselineskip=15pt \small

\noindent{\slshape\bfseries Summary.}
The classical Ray--Knight theorems for Brownian motion determine the law of its local time process either at the first hitting time of a given value $a$ by the local time at the origin, or at the first hitting time of a given position $b$ by Brownian motion. We extend these results by describing the local time process jointly for all $a$ and all $b$, by means of stochastic integral with respect to an appropriate white noise. Our result applies to $\mu$-processes, and has an immediate application: a $\mu$-process is the height process of a Feller continuous-state branching process (CSBP) with immigration (Lambert~\cite{lambert}), whereas a Feller CSBP with immigration satisfies a stochastic differential equation driven by a white noise (Dawson and Li~\cite{dawson-li12}); our result gives an explicit relation between these two descriptions and shows that the stochastic differential equation in question is a reformulation of Tanaka's formula. 
\medskip

\medskip
 

} 


\bigskip

\section{Introduction}

Let $(B_t)_{t\ge 0}$ be standard one-dimensional Brownian motion associated with its completed natural filtration $({\mathscr B}_t)_{t\ge 0}$. Denote by $({\mathfrak L}_t)_{t\ge 0}$ a continuous version of local times of $(B_t)$ at position $0$. Let $\mu \in \r\backslash\{0\} $. The $\mu$-process $X := (X_t)_{t\ge 0}$ is defined as follows:
$$
X_t := |B_t| - \mu {\mathfrak L}_t, \qquad t\ge 0.
$$

\noindent There are two important special cases of $\mu$-processes: Brownian motion ($\mu=1$, this is seen using L\'evy's identity), and the three-dimensional Bessel process ($\mu=-1$, seen by means of L\'evy's and Pitman's identities).

The $\mu$-process, also referred to as perturbed reflecting Brownian motion, has attracted much attention in the nineties: L\'evy's arc sine law, Ray--Knight theorems as well as pathwise uniqueness of doubly perturbed Brownian motion, see for example \cite{legall-yor, yor92, carmona-petit-yor94a, carmona-petit-yor94b, shi-werner95, werner95, chaumont-doney99, davis99, perman96, perman-werner}.

The local time of the $\mu$-process at suitable stopping times, as a process of the space variable, turns out to be a squared Bessel process. This is referred to as a Ray--Knight theorem. More precisely, let us fix $\mu>0$ from now on, so the process $(X_t, \, t\ge 0)$ is recurrent on $\r$.  Since $X$ is a continous semimartingale, we may define   
$$
L(t, \, r) 
:=
\lim_{\varepsilon\to 0} \frac{1}{\varepsilon} \int_0^t  {\bf 1}_{\{ r \le X_s \le r+\varepsilon\} } \, \mathrm{d} s,
\qquad
t\ge 0 , \, r\in \r,
$$

\noindent as the local time of $X$ at time $t$ and position $r$. Moreover, we may and will take a bicontinuous version of local times $L(\cdot, \cdot)$, see \cite{revuz-yor}, Theorem VI.1.7.  Let
\begin{equation}
    \label{def-tau}
    \tau^r_t:=\inf\{ s\ge 0\,:\, L(s,r)> t \} , 
\end{equation}

\noindent be the inverse local time of $X$. Denote by
\begin{equation}
    \label{def-T} 
    T_r :=\inf\{ t\ge 0 : \, X_t = r \},
\end{equation} 

\noindent the hitting time of $r$. The following Ray--Knight theorems were established by Carmona, Petit and Yor \cite{carmona-petit-yor94b} (see also Yor \cite{yor92}, Chapter 9) and by Le Gall and Yor \cite{legall-yor} respectively.
 
\begin{theorem}
 \label{t:RK} 
 
 Fix $\mu>0$. 

(i) (\cite{carmona-petit-yor94b}, \cite{yor92}) Let $a>0$. The process $\big( L(\tau_a^0, -h),\, h\ge 0 \big)$ is a squared Bessel process of dimension $(2-\frac2\mu)$, starting from $a$ and absorbed at $0$. 

(ii) (\cite{legall-yor}) Let $b<0$. The process $\big( L(T_b, b+h),\, 0\le h \le |b| \big)$ is  a squared Bessel process of dimension $\frac2\mu$, starting from $0$ and reflected at $0$.

\end{theorem}

In the special case $\mu=1$: the process $X$ is Brownian motion by L\'evy's identity, so Theorem \ref{t:RK} boils down to the classical Ray--Knight theorem for Brownian motion, originally proved by Ray~\cite{ray} and Knight~\cite{knight} independently.   Werner \cite{werner95} gave an alternative proof of Theorem \ref{t:RK} using a result of Lamperti \cite{lamperti} on semi-stable Markov processes. Perman \cite{perman96} gave another proof of (i) by establishing a path--decomposition result of $X$.   

The aim of this work is to describe the underlying Brownian motion, jointly for all $a$ and $b$, in the local time processes in Theorem \ref{t:RK}.  We do this by means of Tanaka's formula and Walsh's stochastic integral with respect to a white noise $W$; see Theorem \ref{t:main} below. The idea of using Tanaka's formula to prove Ray--Knight theorems is not new, and can be found for example in Jeulin \cite{jeulin85} (for diffusion processes) and in Norris, Rogers and Williams \cite{NRW87} (for Brownian motion with a local time drift); our main contribution is to show how the white noise $W$ explicitly gives the Brownian part jointly for all $a$ and $b$ in Theorem \ref{t:RK}.

The aforementioned white noise $W$ is defined as follows.  For any Borel function $g:\r_+\times \r \to \r$ such that $\int_{\r_+} \d \ell \int_{\r} g^2(\ell,x) \d x<\infty$, let 
\begin{equation}\label{def:W}
W(g):= \int_0^\infty g\big(L(t, X_t), X_t\big)  {\rm sgn}(B_t) \d B_t.
\end{equation}

\noindent It is easily seen that $W$ is a white noise on $\r_+\times \r$; indeed, by the occupation time formula (Exercise VI.1.15 in \cite{revuz-yor}),
$$
\int_0^\infty g^2\big(L(t, X_t), X_t\big)  \d t
=
\int_\r \d x \int_0^\infty g^2\big(L(t, x), x\big)  \d_t L(t, x) 
=
\int_\r \d x \int_0^\infty g^2(\ell, x) \d \ell.
$$

\noindent The exponential martingale for Brownian motion implies that 
$$ 
\e \Big[{\rm e}^{ W(g) - \frac12 \int_0^\infty g^2(L(t, X_t), X_t )\d t}\Big] =1,
$$

\noindent showing that $W(g)$ is a centered Gaussian random variable with variance $\int_{\r_+} \int_{\r} g^2(\ell,x) \d x \d \ell$.

The main result of this work is the following theorem:

\begin{theorem}\label{t:main} Fix $\mu>0$. 
Let $W$ be the white noise defined via \eqref{def:W}.

(i) Almost surely for all $a>0$,
\begin{equation} \label{eq:main1}
L(\tau_a^0,  - h )
=
a - 2 \int_{-h}^0 W\big([0, L(\tau_a^0, x)], \d x\big) + \left(2-{2\over \mu}\right) h, \qquad h\in [0,|I_{\tau_a^0}|],
\end{equation}
\noindent where $I_t:=\inf_{0\le s \le t} X_s$, $t\ge0$, denotes the infimum process of $X$.

(ii) Almost surely for all $b <0$,
\begin{equation} \label{eq:main2}
L(T_{b}, b+h) 
=
  2 \int_b^{b+h} W\big([0, L(T_b, x)], \d x\big)
+ {2\over \mu} h , \qquad h\in [0,|b|].
\end{equation}

\end{theorem}

The precise meaning of stochastic integrals with respect to $W$ is given in Section \ref{s:walsh}. Indeed, we will show that almost surely  \eqref{eq:main1} holds for any fixed $a>0$, hence for all $a$ belonging to a countable dense set of  $\r_+$.  By using the regularity of local times, we may and will choose a version of stochastic integral such that  \eqref{eq:main1}  holds  simultaneously for all $a>0$. The same remark applies to \eqref{eq:main2} as well as to Theorem \ref{t:mainbis} in Section \ref{s:extension}.

It is not surprising, at least in the case (ii), that the local times of a $\mu$-process can be represented as solution of an SDE driven by a white noise. As a matter of fact, by duality (see \cite{werner95}),  the process $(X_{T_{b}-t}-b,\,t\in [0,T_b])$ has the same law as the process $|B_t|+\mu {\mathfrak L}_t$, stopped when leaving $|b|$ for the last time. On one hand, the process $(|B_t|+\mu {\mathfrak L}_t,\, t\ge 0)$ is the height process of a Feller CSBP with immigration (see \cite{lambert},  remark p. 57 in Section 4).  On the other hand, Bertoin and Le Gall showed in \cite{bertoin-legall00} that general CSBPs are  related  to flows of subordinators, constructed in \cite{bertoin-legall-III} critical CSBPs without Gaussian coefficient as solutions of SDEs driven by compensated Poisson random measures. Dawson and Li \cite{dawson-li12} generalized this SDE to include a Gaussian coefficient and possible immigration. Applied to our setting, it is shown that a Feller CSBP with immigration can be constructed as a solution of \eqref{eq:main2}. Theorem \ref{t:main} connects directly the local times of the $\mu$-process to equation \eqref{eq:main2}, without making use of the framework of CSBPs, and in Section \ref{s:pf_thms} we are going to see Theorem \ref{t:RK} as a consequence of Tanaka's formula for $X$.

The rest of the paper is as follows. In Section \ref{s:exc_filtr}, we follow Walsh \cite{walsh78} by introducing the excursion filtration, then make an enlargement of the filtration \`a la Jeulin \cite{jeulin85}. Section \ref{s:walsh} is devoted to study of the martingale measure associated with the white noise $W$. In particular, stochastic integration with respect to $W$ is defined. Theorems \ref{t:main} and \ref{t:RK} are proved in Section \ref{s:pf_thms}. Sections \ref{s:extension} presents analogous results for the $\mu$-process defined on $\r$.

\section{The excursion filtration}
\label{s:exc_filtr}

We first introduce some notation which will be used throughout the paper. 

\begin{notation}\label{n:<h}
Let $x \in \r$. We define the process $X^{-,x}$ obtained by gluing the excursions of $X$ below $x$ as follows. Let, for $t\ge0$, 
$$
A_t^{-,x} := \int_0^t {\bf 1}_{\{X_s \le  x \}} {\rm d} s,\, 
\qquad  
\alpha_t^{-, x}  := \inf\{u>0,\, A_u^{-,x} > t\},
$$

\noindent with the usual convention $\inf \emptyset:=\infty$.  Define 
$$ 
X^{-,x}_t := X_{\alpha_t^{-,x}}, 
\qquad   
t < A_\infty^{-,x}:=\int_0^\infty {\bf 1}_{\{X_s \le  x\} } {\rm d} s.
$$ 

Similarly, we define $A_t^{+,x}$, $\alpha_t^{+,x}$ and $X^{+,x}$ by replacing $X_s \le x$ by $X_s > x$. When the process is denoted by $X$ with some superscript, the analogous quantities  keep the same superscript. For example,  $L^{+,x}(t, y )$ denotes the local time of $X^{+,x}$ at position $y$ and time $t$, and $I_t^{+,x} =\inf_{0\le s\le t} X_s^{+,x}$.  
\end{notation}



\begin{remark} 
Let $x\in \r$. One can reconstruct $X$ from $X^{-,x}$ and $X^{+,x}$ by gluing the excursions of $X^{-,x}$ and of $X^{+,x}$, indexed by their local time. 
\end{remark}

The following proposition  is adapted from Section 8.5 of \cite{yor92}.

\begin{proposition}\label{p:decomp}
Let $x\le 0$. 

(i) Define the filtration $({\mathscr F}^{+,x}_u)_{u\ge 0}$ by ${\mathscr F}_u^{+,x} := \sigma(X_s^{+,x},\, s \in [0,u])$ and the process
$$
\beta_u^{+,x} := \int_0^{\alpha_u^{+,x}} {\bf 1}_{\{X_s > x\}} {\rm sgn}(B_s) \d B_s, \qquad u\ge0.
$$

\noindent Then $\beta^{+,x}$ is $({\mathscr F}^{+,x}_u)$-Brownian motion and $X^{+,x}$ is an $({\mathscr F}^{+,x}_u)$-semimartingale with decomposition
\begin{equation}\label{eq:decomp+}
X_u^{+,x} =  \beta_u^{+,x} - {1-\mu \over \mu}I^{+,x}_u + {1\over 2} L^{+,x}(u,x),
\qquad
u\ge 0.
\end{equation}

(ii) Define the filtration $({\mathscr F}^{-,x}_u)_{u\ge 0}$ by ${\mathscr F}_u^{-,x} := \sigma( X_s^{-,x},\, s \in [0,u])$ and the process
$$
\beta_u^{-,x} := \int_{0}^{\alpha_u^{-,x}} {\bf 1}_{\{X_s \le x\}} {\rm sgn}(B_s) \d B_s,  \qquad u\ge0.
$$

\noindent Then $\beta^{-,x}$ is $({\mathscr F}^{-,x}_u)$-Brownian motion and $X^{-,x}$ is an $({\mathscr F}^{-,x}_u)$-semimartingale with decomposition
\begin{equation}\label{eq:decomp-}
X_u^{-,x} = x + \beta_u^{-,x} - {1-\mu \over \mu} (I^{-,x}_u-x) -{1\over 2} L^{-,x}(u,x),
\qquad
u\ge 0.
\end{equation}

(iii) The Brownian motions $\beta^{+,x}$ and $\beta^{-,x}$ are independent.
\end{proposition}

{\noindent\it Proof}. By Tanaka's formula, 
$$
(X_t-x)^+= (X_0-x)^+ + \int_0^t {\bf 1}_{\{X_s > x\}} \d X_s + \frac12 L(t, x).
$$

\noindent Take $t = \alpha^{+,x}_u$. We get
\begin{equation}\label{eq:tanaka1}
X_u^{+,x}=   \int_0^{\alpha^{+,x}_u} {\bf 1}_{\{X_s > x\}} \d X_s 
+ \frac12 L(\alpha^{+,x}_u, x).
\end{equation}

\noindent Moreover, $\d X_s = \d |B_s| - \mu \d {\mathfrak L}_s = {\rm sgn}(B_s)\d B_s + (1-\mu) \d {\mathfrak L}_s$ by another application of Tanaka's formula. Also observe that $I_t=-\mu {\mathfrak L}_t$, hence \begin{equation}\label{dXs}\d X_s = {\rm sgn}(B_s)\d B_s - {1-\mu\over \mu} \d I_s, \qquad s\ge0. \end{equation}

\noindent Therefore, 
$$
X_u^{+,x}=   \int_0^{\alpha^{+,x}_u} {\bf 1}_{\{X_s > x\}} {\rm sgn}(B_s)\d B_s  - {1-\mu\over \mu}\int_0^{\alpha^{+,x}_u}{\bf 1}_{\{X_s > x\}} \d I_s
+ \frac12 L(\alpha^{+,x}_u, x).
$$

\noindent We  notice    that $L^{+,x}(u,r) = L(\alpha_u^{+,x},r)$ for any $r\in [x, \infty)$ and $u\ge 0$. On the other hand, $
\int_0^{\alpha^{+,x}_u}{\bf 1}_{\{X_s > x\}} \d I_s=I_{\alpha^{+,x}_u\land T_x} $ which is also the infimum of $X^{+,x}$ on the time interval $[0,u]$. This yields \eqref{eq:decomp+}. This equation also implies that $\beta^{+,x}$ is adapted to ${\mathscr F}^{+,x}$. Moreover,  from the definition of  $\beta^{+,x}_u$ and Proposition V.1.5 of \cite{revuz-yor}, $\beta^{+,x}_u$ is a $({\mathscr B}_{\alpha_u^{+,x}})$-continuous martingale with $\langle \beta^{+,x}, \beta^{+,x}\rangle_u=u$, hence $({\mathscr B}_{\alpha_u^{+,x}})$-Brownian motion. Since ${\mathscr F}_u^{+,x} \subset {\mathscr B}_{\alpha_u^{+,x}}$, we deduce that $\beta^{+,x}_u$ is also $({\mathscr F}_u^{+,x})$-Brownian motion. 
This proves (i).

The proof of (ii) is similar.  Tanaka's formula applied to $(X_{t+T_x}-x)^-$ with $t=\alpha^{-,x}_u - T_x$ implies that
$$
X_u^{-,x}= x +  \int_{T_x}^{\alpha^{-,x}_u} {\bf 1}_{\{X_s \le x\}} {\rm sgn}(B_s)\d B_s  - {1-\mu\over \mu}\int_{T_x}^{\alpha^{-,x}_u}{\bf 1}_{\{X_s \le x\}} \d I_s
- \frac12 L(\alpha^{-,x}_u, x).
$$

\noindent We  observe that $L(\alpha^{-,x}_u, x)=L^{-,x}(u,x)$ and $\int_{T_x}^{\alpha^{-,x}_u}{\bf 1}_{\{X_s \le x\}} \d I_s = I_{\alpha^{-,x}_u}-x$ while $I_{\alpha^{-,x}_u}= I^{-,x}_u$ which gives \eqref{eq:decomp-}. We conclude as for (i). The statement (iii) is a consequence of Knight's theorem on orthogonal martingales. $\Box$

\bigskip

The following result is well-known. It has been proved in Section 8.5 of \cite{yor92} when $\mu\in(0,2)$, in \cite{werner95} and in \cite{perman-werner}. Here, following \cite{yor92}, we choose to see it as a consequence of Proposition \ref{p:decomp}.

\begin{corollary}\label{c:independence}
Let $x\le 0$. The processes $X^{+,x}$ and $X^{-,x}$ are independent.
\end{corollary}

{\noindent\it Proof}. By Proposition \ref{p:decomp}~(iii), the martingale parts of $X^{+,x}$ and $X^{-,x}$, namely $\beta^{+,x}$ and $\beta^{-,x}$, are independent. It remains to see that $X^{+,x}$ is measurable with respect to $\beta^{+,x}$ and $X^{-,x}$ with respect to $\beta^{-,x}$, which was established by Chaumont and Doney \cite{chaumont-doney99} and Davis \cite{davis99}. $\Box$

\bigskip

The excursion filtration, introduced by Walsh \cite{walsh78}, is defined as
$$
{\mathcal E}^+_x := {\mathscr F}_\infty^{+,x} = \sigma(X_s^{+,x},\, s \ge0) ,\qquad x\in \r.
$$

\noindent Similarly we define ${\mathcal E}^-_x:={\mathscr F}_\infty^{-,x}=\sigma(X_s^{-,x},\, s \ge0)$ for $x\in \r$.  It is routine to check, using the time-changes $\alpha^{-,x}$ and $\alpha^{+,x}$, that  ${\mathcal E}^-_x$ is increasing in $x$ whereas ${\mathcal E}^+_x$ is decreasing. \footnote{We will be implicitly working with a right-continuous (and complete) version of the filtrations $({\mathcal E}^{-}_x)_{x\in \r}$ and $({\mathcal E}^{+}_x)_{x\in \r}$ --- if necessary, by means of the procedure of usual augmention, as described in Section I.4 of Revuz and Yor~\cite{revuz-yor}. In our work, we only study {\it continuous} martingales, which are also martingales with respect to augmented filtrations. The same remark applies to the $\mu$-process defined on $\r$ in Section \ref{s:extension}.}

Define, for $u\ge 0$,
$$
{\mathscr G}^{+,x}_u:= \sigma( {\mathscr F}^{+,x}_u, {\mathcal E}^-_x),
\qquad
{\mathscr G}^{-,x}_u:= \sigma( {\mathscr F}^{-,x}_u, {\mathcal E}^+_x) .
$$

\noindent The idea of such an enlargement of filtrations goes back at least to Jeulin \cite{jeulin85}.

\begin{corollary}\label{c:sgnB}
 Consider a   random function $g(\ell, y)=g(\ell, y,\omega)$ such that the process $t \mapsto g(L(t, X_t), X_t)$ is $({\mathscr B}_t)$-progressively measurable and $\e\big[\int_{\r_+} \d \ell \int_\r g(\ell, y)^2 \d y \big] <\infty$. Fix $x \le 0$.

(i)
 The process $u \mapsto g\big(L^{+,x}(u,X^{+,x}_u),X^{+,x}_u\big)$ is $({\mathscr G}^{+,x}_u)$-progressive  and almost surely,
\begin{equation}\label{eq:sgnB+}
 \int_0^\infty g\big(L(t, X_t), X_t\big) {\bf 1}_{\{X_t > x \}}   {\rm sgn}(B_t) \d B_t
 =
 \int_0^\infty g\big(L^{+, x} (u,X_u^{+, x}), X_u^{+,x}\big) {\rm d} \beta_u^{+,x}. 
\end{equation}

(ii) The process $u \mapsto g\big(L^{-,x}(u,X^{-,x}_u),X^{-,x}_u\big)$ is $({\mathscr G}^{-,x}_u)$-progressive  and almost surely,
\begin{equation}\label{eq:sgnB-}
\int_0^\infty g\big(L(t, X_t), X_t\big) {\bf 1}_{\{X_t \le x \}}   {\rm sgn}(B_t) \d B_t
 =
 \int_0^\infty g\big(L^{-, x} (u,X_u^{-, x}), X_u^{-,x} \big) {\rm d} \beta_u^{-,x}. 
\end{equation} 
 \end{corollary}

{\noindent\it Proof}. We prove (i). The process $u \mapsto g(L^{+,x}(u,X^{+,x}_u),X^{+,x}_u)$ is  $({\mathscr B}_{\alpha^{+,x}_u})$-progressive (Proposition V.1.4,  \cite{revuz-yor}). Therefore it is also progressive with respect to $({\mathscr G}^{+,x}_u)$ because the latter filtration is larger. We prove now \eqref{eq:sgnB+}. By a time-change (Proposition V.1.4, \cite{revuz-yor}), 
$$
\int_{0}^{\alpha_t^{+,x}} g\big(L(s,X_s) ,X_s\big){\bf 1}_{\{ X_s > x\}} {\rm sgn}(B_s)\d B_s=\int_0^t g\big(L^{+,x} (u,X_u^{+,x}) , X_u^{+,x}  \big)  \d \beta_u^{+,x}.
$$

\noindent Letting $t\to \infty$ yields \eqref{eq:sgnB+}. Statement (ii) is proved similarly.   $\Box$

\section{The martingale measure associated with $W$}
\label{s:walsh}

 Recall the definition of the white noise $W$ in \eqref{def:W}. For any  Borel set $A$ of $\r_+$ with finite Lebesgue measure  and $r\ge 0$, we define   \begin{equation}\label{def:M}
M_r(A):=W({\bf 1}_{A\times [-r,0]}).
\end{equation}

\begin{proposition}\label{p:martingale}
In the setting of Walsh \cite{walsh84}, $(M_r, \, r\ge 0 )$ is a continuous martingale measure with respect to the filtration $({\mathcal E}^+_{-r},\,r\ge 0)$. 
\end{proposition}

{\noindent\it Proof}. Since $W$ is a white noise, it suffices to show that $M_r$ is measurable with respect to ${\mathcal E}_{-r}^+$ and that $M_{s}-M_r$ is independent of ${\mathcal E}_{-r}^+$ for any $0\le r < s$. The first statement comes from \eqref{eq:sgnB+} applied to $x=-r$ and $g(\ell,y)={\bf 1}_{A\times [-r,0]}(\ell,y)$ for a Borel set $A$ with finite Lebesgue measure. The second statement comes from \eqref{eq:sgnB-} applied to $x=-r$ and $g(\ell,y)={\bf 1}_{A \times [-s,-r)}(\ell,y)$.  Since the processes $X^{-,-r}$ and $\beta^{-, -r}$ are independent of ${\mathcal E}_{-r}^+$,  the proposition follows. $\Box$

\bigskip


We are going to extend \eqref{def:W}, seen as an equality for deterministic functions $g$,  to random functions. To this end, we first recall the construction by Walsh in \cite{walsh84} of stochastic integral with respect to the martingale measure $M$. A (random) function $f$ is said to be elementary if it is of the form $f( \ell , x ):= Z \, {\bf 1}_{[a, b)}(x) \, {\bf 1}_A(\ell)$, where $a <b\le0$,  $A\subset \r_+$ is a Borel set of finite Lebesgue measure,  and $Z$ is a bounded     ${\cal E}^{+}_{b}$-measurable real-valued random variable.  Denote by $f\cdot M$ the stochastic integral with respect to $M$: 
$$
f\cdot M:= Z\, (M_{|a|}(A)- M_{|b|}(A))= Z \, W({\bf 1}_{A \times [a, b)}).
$$

A simple function is a (finite) linear combination   of elementary functions. We extend by linearity the definition of $f\cdot M$ to   simple functions $f$ and furthermore by   isometry to any  $f \in {\mathscr L}^2$, where  ${\mathscr L}^2$ denotes  the space of $({\mathcal E}^+_{-r},\,r\ge 0)$-predictable  and square-integrable functions, defined as the closure of the space of simple functions under the  norm: 
$$
\|f\| := \Big[ \e\Big( \int_{\r_+} \d \ell \int_{\r_-}  f^2(\ell, x) \d x \Big) \Big]^{1/2} .
$$
 
For any $f \in {\mathscr L}^2$, $f\cdot M$ is a centered random variable with $\e[(f\cdot M)^2] = \|f\|^2$. We write $f\cdot M\equiv \int_{\r_+\times \r_-} f(\ell, x) W(\d \ell, \d x)$ and for any $r\ge0$, $(f{\bf 1}_{\r_+\times [-r, 0]}) \cdot M\equiv \int_{\r_+\times [-r, 0]} f(\ell, x) W(\d \ell, \d x)$. The latter, if furthermore $f$ is of form $f(\ell, x)= {\bf 1}_{\{0\le \ell \le \sigma_x\}} \eta_x$, will be re-written as $\int_{-r}^0 \eta_x W\big([0, \sigma_x], \d x\big)$. By the construction of stochastic integral and Proposition \ref{p:martingale}, $\int_{\r_+\times [-r, 0]} f(\ell, x) W(\d \ell, \d x)$  is a continuous martingale with respect to the filtration $({\mathcal E}^+_{-r},\,r\ge 0)$, of quadratic variation process $\int_{\r_+\times [-r, 0]} f^2(\ell, x)  \d \ell  \d x$.

\begin{proposition}\label{p:main}
Take $g \in {\mathscr L}^2$ such that $s\mapsto g\big(L(s, X_s), X_s\big)$ admits a version which is progressive with respect to the Brownian filtration $({\mathscr B}_s)$. \footnote{By {\it version} we mean  a $({\mathscr B}_s)$-progressive process $(h_s)$ such that  $\int_0^\infty {\bf 1}_{\{g(L(s, X_s), X_s ) \neq h_s\}} \d s =0$ a.s.}  Then
$$
g\cdot M 
=
\int_0^\infty g\big(L(s, X_s), X_s\big)  {\rm sgn}(B_s) \d B_s 
\qquad\hbox{\rm a.s.}
$$
\end{proposition}
 
 {\noindent\it Proof}. By definition of ${\mathscr L}^2$, there exists a sequence of simple functions $g_n$ such that $\| g- g_n\|\to 0$ as $n\to\infty$. By isometry, $ g \cdot M - g_n\cdot M \to 0 $ in $L^2$.  Since $g_n$ is a simple function, $g_n$ is of the form $$ g_n(\ell, x)= \sum_{k, j=1}^\infty  Z^n_{k,j} {\bf 1}_{[a^n_k, a^n_{k-1})}(x)  {\bf 1}_{A^n_j}(\ell),$$

\noindent where for each $n$, $0=a^n_0> ...>a^n_k>a^n_{k+1}>... $ is a decreasing   sequence such that $a^n_k\to -\infty$ as $k\to \infty$,  $(A^n_j)_{j\ge 1}$ is  
a collection of (nonrandom) pairwise disjoint Borel subsets of $\r_+$ with finite Lebesgue measures, and for any $k, j\ge 1$, $Z^n_{k,j}$ is a bounded ${\cal E}^+_{a^n_{k-1}}$-measurable random variable.  Moreover for all large $k, j$, $Z^n_{k,j}=0$, which means the above double sum  runs in fact over a finite index set  of $k$ and $j$. 

Note that for a.e.\ $z\le 0$, $g(\cdot, z)$ is measurable with respect to ${\cal E}^+_z$ (as $g_n$ is). We may (and will) take a version of $g$ such that $g(\cdot, z)$ is measurable with respect to ${\cal E}^+_z$ for all $z \le 0$.

By applying \eqref{def:W} and \eqref{def:M}, we deduce from the linearity of the integral  that 
$$
g_n \cdot M
=
\sum_{k, j=1}^\infty  Z^n_{k,j}  \int_0^\infty {\bf 1}_{[a^n_k, a^n_{k-1})}(X_s)  {\bf 1}_{A^n_j}(L(s,X_s)) {\rm sgn}(B_s) \d B_s.$$

\noindent Note that we can (and we will) take $(a^n_k,a^n_{k-1})$ instead of $[a^n_k, a^n_{k-1})$ without changing the value of $g_n \cdot M$.  By \eqref{eq:sgnB+} and \eqref{eq:sgnB-} respectively, 
\begin{eqnarray}\nonumber
&& \int_0^\infty {\bf 1}_{(a^n_k, a^n_{k-1})}(X_s)  {\bf 1}_{A^n_j}(L(s,X_s)) {\rm sgn}(B_s) \d B_s
\\ 
&=& 
\int_0^\infty {\bf 1}_{(a^n_k, a^n_{k-1})}(X_u^{-,x}) {\bf 1}_{A^n_j}(L^{-,x}(u,X_u^{-,x})) \d \beta_u^{-,x} \label{eq:gn-}
\\
&=&
\int_0^\infty {\bf 1}_{(a^n_k, a^n_{k-1})}(X_u^{+,y}) {\bf 1}_{A^n_j}(L^{+,y}(u,X_u^{+,y})) \d \beta_u^{+,y}  \label{eq:gn+}
\end{eqnarray}

\noindent with $x=a_{k-1}^n$ and $y=a_k^n$. Similarly,
\begin{eqnarray}\nonumber
&& \int_0^\infty g\big(L(s, X_s), X_s\big) {\bf 1}_{(a^n_k, a^n_{k-1})}(X_s) {\rm sgn}(B_s) \d B_s \\ 
&=&
\int_0^\infty g\big(L^{-,x}(u,X_u^{-,x}),X_u^{-,x}\big){\bf 1}_{(a^n_k, a^n_{k-1})}(X_u^{-,x}) \d \beta_u^{-,x},\label{eq:g-}  \\
&=& 
\int_0^\infty g\big(L^{+,y}(u,X_u^{+,y}),X_u^{+,y}\big){\bf 1}_{(a^n_k, a^n_{k-1})}(X_u^{+,y}) \d \beta_u^{+,y} \label{eq:g+}
\end{eqnarray} 

\noindent with $x=a_{k-1}^n$ and $y=a_k^n$ as in    \eqref{eq:gn-} and \eqref{eq:gn+}.  Write $$I(g):= \int_0^\infty g\big(L(s, X_s), X_s\big)  {\rm sgn}(B_s) \d B_s. $$

\noindent Then  $$g_n\cdot M - I(g)= \sum_{k=1}^\infty \Delta_n(k),$$

\noindent where, from \eqref{eq:gn+} and \eqref{eq:g+}, 
\begin{eqnarray*} \Delta_n(k)
&:=& 
 \sum_{j=1}^\infty   Z^n_{k,j} \int_0^\infty {\bf 1}_{(a^n_k, a^n_{k-1})}(X_u^{+,y}) {\bf 1}_{A^n_j}(L^{+,y}(u,X_u^{+,y})) \d \beta_u^{+,y}   
  \\
  && - \int_0^\infty g\big(L^{+,y}(u,X_u^{+,y}),X_u^{+,y}\big){\bf 1}_{(a^n_k, a^n_{k-1})}(X_u^{+,y}) \d \beta_u^{+,y},
\end{eqnarray*}

\noindent where as before $y=a_k^n$. Since $Z^n_{k,j}$ is ${\cal E}^+_{a^n_{k-1}}$-measurable, hence ${\cal E}^+_{a^n_k}$-measurable, we see that the sum over $j$ in the definition of $\Delta_n(k)$ is ${\cal E}^+_{a^n_k}$-measurable. For the last integral in $\Delta_n(k)$, we use the fact that for $z\le 0$, $g(\cdot, z)$ is measurable with respect to ${\cal E}^+_z$.  Note that $X_u^{+,y}$ and $L^{+,y}(\cdot, \cdot)$ are  measurable with respect to ${\mathcal E}^+_y$. Since $X_u^{+,y} \ge y$ and ${\mathcal E}^+_z$ decreases on $z$, we deduce that $g(\cdot,X_u^{+,y})$ is  ${\mathcal E}^+_y$-mesurable, and so is $g\big(L^{+,y}(u,X_u^{+,y}),X_u^{+,y}\big)$. It follows that $\Delta_n(k)$ is measurable with respect to ${\mathcal E}^+_{a_k^n}$. 

Now we prove that $(\Delta_n(k))_{k\ge1}$ is a martingale difference sequence with respect to the filtration  $({\mathcal E}^+_{a_k^n})_{k\ge1}$. Indeed, using  \eqref{eq:gn-} and \eqref{eq:g-} instead of \eqref{eq:gn+} and \eqref{eq:g+}, one can also write
\begin{eqnarray*} \Delta_n(k)
&=& 
 \sum_{j=1}^\infty   Z^n_{k,j} \int_0^\infty {\bf 1}_{(a^n_k, a^n_{k-1})}(X_u^{-,x}) {\bf 1}_{A^n_j}(L^{-,x}(u,X_u^{-,x})) \d \beta_u^{-,x}   
  \\
  && - \int_0^\infty g\big(L^{-,x}(u,X_u^{-,x}),X_u^{-,x}\big){\bf 1}_{(a^n_k, a^n_{k-1})}(X_u^{-,x}) \d \beta_u^{-,x},
\end{eqnarray*}

\noindent with  $x=a_{k-1}^n$. Recall from Proposition \ref{p:decomp} that $\beta^{-,x}$ is $({\mathscr F}^{-,x}_u)$-Brownian motion, which is independent of ${\mathcal E}^+_x$ by Corollary \ref{c:independence}. Then $\beta^{-,x}$ can be seen as $({\mathscr G}^{-,x}_u)$-Brownian motion. By Corollary \ref{c:sgnB}~(ii), $g\big(L^{-,x}(u,X_u^{-,x}),X_u^{-,x}\big)$ is progressive with respect to the filtration $({\mathscr G}^{-,x}_u)$, while 
$$
{\bf 1}_{(a^n_k, a^n_{k-1})}(X_u^{-,x}) \,  g_n  \big(L^{-,x}(u,X_u^{-,x}),X_u^{-,x}\big) 
=
\sum_{j=1}^\infty   Z^n_{k,j} {\bf 1}_{(a^n_k, a^n_{k-1})}(X_u^{-,x}) {\bf 1}_{A^n_j}(L^{-,x}(u,X_u^{-,x}))
$$

\noindent is $({\mathscr G}^{-,x}_u)$-progressive as well. Therefore, one can write 
\begin{equation}\label{eq:Delta}
\Delta_n(k)
=
\int_0^\infty {\bf 1}_{(a^n_k, a^n_{k-1})}(X_u^{-,x}) \, (g_n -g)\big(L^{-,x}(u,X_u^{-,x}),X_u^{-,x}\big) \d \beta_u^{-,x},
\end{equation}

\noindent with $x=a_{k-1}^n$.  It follows that (since ${\mathscr G}_0^{-,x}={\mathcal E}^+_x$)
$$
\e[ \Delta_n(k) \mid {\mathcal E}^+_{a^n_{k-1}} ]=0.
$$

\noindent In other words, the process $j\to \sum_{k=1}^j \Delta_n(k)$ is a martingale and we have
$$
\e[(g_n\cdot M - I(g))^2]= \sum_{k=1}^\infty \e[\Delta_n(k)^2].
$$

\noindent From \eqref{eq:Delta}, we get 
\begin{eqnarray*}
\e[\Delta_n(k)^2]
&=&
 \e\left[ \int_0^\infty  {\bf 1}_{(a^n_k, a^n_{k-1})}(X_u^{-,x}) \, \big((g_n -g)(L^{-,x}(u,X_u^{-,x}),X_u^{-,x})\big)^2 \d u \right]
\\
&=&  \e\left[\int_{a^n_k}^{a^n_{k-1}} (g_n(\ell, z) -g(\ell,z))^2 \d \ell \d z \right],
\end{eqnarray*}

\noindent where in the second equality we have used the occupation time formula. It follows that 
$$
\e[(g_n\cdot M - I(g))^2]= \|g_n-g\|^2
$$

\noindent  which goes to $0$ as $n\to\infty$. Thus we get that $g\cdot M = I(g)$. $\Box$

\begin{remark}\label{r:E-}  Fix $b<0$. Similarly to \eqref{def:M} we may define a martingale measure $\widehat M$ by  $$
\widehat  M_r(A):=W({\bf 1}_{A\times [b, b+r]}) , $$

\noindent for any Borel set $A\subset \r_+$ of finite Lebesgue measure  and $r\ge0$.  The analogs of  Propositions \ref{p:martingale} and \ref{p:main} hold for $\widehat M$. Specifically, $(\widehat  M_r, \, r\ge 0 )$ is a martingale measure with respect to the filtration $({\mathcal E}^-_{b+r},\,r\ge 0)$. Moreover, we may define in a similar way the stochastic integral $f \cdot \widehat M$  for any $f \in \widehat  {\mathscr L}^2$, where  $\widehat {\mathscr L}^2$ denotes the space of $({\mathcal E}^-_{b+r})_{r\ge0}$-predictable functions $f$  such that $\e\big[ \int_{\r_+} \d \ell \int_{[b, \infty)} f^2(\ell, x) \d x   \big] < \infty$.  Then for any $g \in \widehat  {\mathscr L}^2$ such that   $s\mapsto g(L(s, X_s), X_s)$ admits a version which is  progressive with respect to the   Brownian filtration $({\mathscr B}_s)$, we have
\begin{equation}\label{r:eq} 
g\cdot \widehat  M 
=\int_0^\infty g\big(L(s, X_s), X_s\big)   \d X_s - \frac{1-\mu}{\mu} \int_{-\infty}^0 g(0, x) \d x,
\qquad\hbox{\rm a.s.}  
\end{equation}

With a slight abuse of notation, we shall write $g\cdot \widehat M\equiv \int_{\r_+ \times [b, \infty)} g(\ell, x) W(\d \ell, \d x) $ and $(g{\bf 1}_{\r_+\times [b, t)}) \cdot \widehat M\equiv \int_{\r_+ \times [b, t)}  g(  \ell, x) W(\d \ell, \d x) $ for any $t\ge b$. Then $\int_{\r_+ \times [b, b+r)}  g(  \ell, x) W(\d \ell, \d x)$, is an $({\mathcal E}^-_{b+r})$-continuous martingale with quadratic variation process $\int_{\r_+ \times [b, b+r)}  g^2(\ell, x)  \d \ell  \d x$ for $r\ge0$. 

\end{remark}

\section{Proofs of Theorems \ref{t:main} and \ref{t:RK}}
\label{s:pf_thms} 
 
By Tanaka's formula, for any $r \ge0$ and $x \in \r$,
\begin{equation}\label{eq:tanaka}
(X_r-x)^-=(-x)^- - \int_0^r {\bf 1}_{\{ X_s \le x\}} \d X_s + \frac12 L(r,x).
\end{equation}

{\noindent\it Part (i) of Theorems \ref{t:main} and \ref{t:RK}:}  Applying \eqref{eq:tanaka} to $r= \tau_a^0$ gives that $L(\tau_a^0, x)= 2 \int_0^{\tau_a^0} {\bf 1}_{\{ X_s \le x\}} \d X_s$ for all $x\in \r$. Let $h\ge0$. Taking $x=-h$ and $x=0$,  and using the fact that $L(\tau_a^0, 0)=a$,  we obtain that  
\begin{eqnarray}
L(\tau_a^0,  -h)
&=& a -  2 \int_0^{\tau_a^0} {\bf 1}_{\{ -h < X_s \le  0\}} \d X_s  \nonumber
\\
&=&
a -  2 \int_0^{\tau_a^0} {\bf 1}_{\{ -h  < X_s \le  0\}} {\rm sgn}(B_s)\d B_s +  \Big(2- \frac{2}{\mu}\Big) \min(h , |I_{\tau_a^0}|),  \label{eq:tanaka-tau}
\end{eqnarray}

\noindent where the second equality follows from \eqref{dXs} and the fact that $\int_0^{\tau_a^0} {\bf 1}_{\{ -h  < X_s \le  0\}}   \d I_s =\int_0^{\tau_a^0} {\bf 1}_{\{I_s>-h\}}   \d I_s= I_{\min(T_{-h}, I_{\tau_a^0})}= - \min(h , |I_{\tau_a^0}|)$.

 To deal with the stochastic integral with respect to $(B_s)$ in \eqref{eq:tanaka-tau}, we shall use Proposition \ref{p:main}. First we remark  that for $x\le 0$, $L(\tau_a^0, x)$ is measurable with respect to  ${\mathcal E}^+_x$. In fact,  let $u:= A^{+, x}_{\tau^0_a}$. Then $\alpha^{+,x}_u=\tau_a^0$ as  $\tau_a^0$ is an increasing time for $A^{+,x}$. Therefore $L(\tau^0_a,x)=L^{+,x}(A^{+,x}_{\tau^0_a},x)$. Since $\{A^{+,x}_{\tau^0_a}\ge t\}= \{L(\alpha^{+,x}_t, 0) \le a\}= \{L^{+,x}(t, 0) \le a\}$, we obtain that $L(\tau^0_a, x)$  is measurable with respect to ${\mathcal E}^+_x$. 
 
 Let $g(\ell, x):= {\bf 1}_{  \{0 \le \ell \le L(\tau_a^0, x)\} \cap \{-h < x \le 0\}}$.  Using the continuity of local times $L(\tau_a^0, x)$ on $x$ and the   fact that $\e[\int_{\r_+\times \r} g^2(\ell, x) \d \ell \d x]< \infty$,  we get that   $g \in {\mathscr L}^2$.

Observe that a.s., $g(L(s, X_s), X_s)  = {\bf 1}_{\{0\le s \le \tau_a^0, \, -h< X_s\le 0\}}$, $\d s$-a.e. This follows from the fact that $\int_0^\infty {\bf 1}_{\{L(s, X_s) \le L(\tau_a^0, X_s), s > \tau_a^0\}  } \d s= \int_{-\infty}^\infty \d x \int_{\tau_a^0}^\infty {\bf 1}_{\{L(s, x) \le L(\tau_a^0, x)\}} \d_s L(s, x)=0$, by the occupation time formula.  
 Then $s\mapsto g(L(s, X_s), X_s)$ admits a version which is $({\mathscr B}_s)$-progressive and we are entitled to apply Proposition \ref{p:main} to see that 
$$
\int_0^{\tau_a^0} {\bf 1}_{\{ -h  < X_s \le  0\}} {\rm sgn}(B_s)\d B_s
=
\int_{\r_+\times (-h,0] } {\bf 1}_{\{ 0 \le \ell \le L(\tau_a^0, x) \}} W(\d \ell, \d x)
=
\int_{-h}^0  W\big([0, L(\tau_a^0, x)], \d x\big).
$$

\noindent  In view of \eqref{eq:tanaka-tau}, this yields Theorem \ref{t:main}~(i) for each fixed $a>0$ and $h\in [0, |I_{\tau_a^0}|]$. Since the processes are continuous in $h$ and c\`adl\`ag in $a$, they coincide except on a null set. 

Theorem \ref{t:RK}~(i) follows quite simply from Theorem \ref{t:main}~(i): for given $a>0$, since $h\mapsto \int_{-h}^0 W\big([0, L(\tau_a^0, x)], \d x\big)$ is an $({\mathcal E}^+_{-h})_{h\ge0}$-continuous martingale with quadratic variation process $\int_{-h}^0 L(\tau_a^0, x) \d x$, it follows from the Dambis--Dubins--Schwarz theorem that there exists $({\mathcal E}^+_{-h})$-Brownian motion $\gamma$ such that $\int_{-h}^0  W\big([0, L(\tau_a^0, x)], \d x\big) = \int_0^h \sqrt{ L(\tau_a^0, s)} \d \gamma_s$. Going back to \eqref{eq:tanaka-tau}, we see that for all $0\le h \le \inf\{s\ge0: L(\tau_a^0, -s) =0\}= |I_{\tau_a^0}|$, 
$$
L(\tau_a^0,  - h )
=
a - 2 \int_0^h \sqrt{ L(\tau_a^0, s)} \d \gamma_s + \Big(2- \frac{2}{\mu}\Big) h,
$$

\noindent proving Theorem \ref{t:RK}~(i). 
 
\medskip
{\noindent\it Part (ii) of Theorems \ref{t:main} and \ref{t:RK}:}  Let $b<0$ and $T_b:=\inf\{t\ge0: X_t=b\}$. For $h \in [0,|b|]$, we get from \eqref{eq:tanaka} that
\begin{eqnarray}
L(T_b, b+h) 
&=&
2h + 2 \int_0^{T_b} {\bf 1}_{\{X_s \le b+h\}} d X_s \nonumber
\\
&=&
  2 \int_0^{T_b} {\bf 1}_{\{X_s \le b+h\}} {\rm sgn}(B_s)\d B_s  + {2\over \mu} h , 
\label{eq:tanaka-tb}
\end{eqnarray}

\noindent where the second equality follows from \eqref{dXs} again  and the fact that $\int_0^{T_b} {\bf 1}_{\{X_s \le b+h\}}   \d I_s =\int_0^{T_b} {\bf 1}_{\{I_s \le b+h\}}   \d I_s = -h$.

The main difference with  Part (i)  is the measurability. As a matter of fact, for any $x\ge b$, $L(T_b, x)$ is ${\mathcal E}^-_x$-measurable:  observe that $L(T_b, x)= L^{-,x}(A^{-, x}_{T_b}, x)$ and for any $t\ge0$, $\{A^{-, x}_{T_b}> t\}= \{T_b > \alpha^{-,x}_t\}=\{\inf_{0\le s \le \alpha^{-,x}_t}X_s > b\}= \{\inf_{0\le s \le t} X^{-,x}_s > b\}$ is ${\mathcal E}^-_x$-mesurable.

  Let $g(\ell, x):= {\bf 1}_{ \{0< \ell \le L(T_b,x)\} \times \{b \le x \le b+h\} }$. We can check as in (i) that we may apply \eqref{r:eq} to get that 
$$
\int_0^{T_b} {\bf 1}_{\{X_s \le b+h\}} {\rm sgn}(B_s)\d B_s
=
 \int_{ \r_+ \times [b,b+h]} {\bf 1}_{\{ 0\le \ell \le  L(T_b, b+h) \}}W(\d \ell, \d x)
 =
 \int_b^{b+h}  W\big([0, L(T_b, x)], \d x\big),
$$

\noindent proving, in view of \eqref{eq:tanaka-tb},    Theorem \ref{t:main} (ii). Furthermore by Remark \ref{r:E-}, the process $h \to \int_b^{b+h}  W\big([0, L(T_b, x)], \d x\big)$ is an $({\mathcal E}^-_{b+h})$-continuous martingale with quadratic variation process $  \int_{b}^{b+h} L(T_b, x) \d x$. This easily yields   Theorem \ref{t:RK} (ii). $\Box$

\section{Extension to the two-sided $\mu$-process}\label{s:extension}

In this Section, we shall explore the strong Markov property at the hitting times of a $\mu$-process defined on $\r$ and  present an analogue  of Theorem \ref{t:main}. This result, apart from its own interest,  will be useful in a forthcoming work on the duality of Jacobi stochastic flows.
 
 Let $(B_t)_{t\in\r}$ be a two-sided Brownian motion, which means that for $t\le 0$, $B_t=B'_{-t}$, where $B'$ is a standard Brownian motion independent of $(B_t)_{t\ge0}$.  Denote by $({\mathfrak L}'_t)_{t\ge0}$ the local time process at position zero of $B'$.

  Recall that $X_t=|B_t| - \mu {\mathfrak L}_t$, for $t\ge 0$.  For $t\le 0$, we let $X_t:=|B_{-t}'|+\mu {\mathfrak L}'_{-t}$. We call $(X_t)_{t\in \r}$ a two-sided $\mu$-process. Fix $\mu >0$.  Notice that $X_t\to \infty$ as $t \to -\infty$, and $T_r<0$ when $r>0$.

   We naturally extend the notation $T_r:=\inf\{t\in\r\,:\, X_t=r\}$ for  $r \in \r$, $$L(t,x):= \lim_{\varepsilon \to 0} \frac1\varepsilon \int_{-\infty}^t 1_{\{ x \le X_s \le x+\varepsilon\}} \d s, \qquad x \in \r,  $$ the local time accumulated by $(X_t,\, t\in \r)$ at position $x$ up to time $t$, and $$\tau_a^x:= \inf\{t \in \r: L(t, x) > a\}, \qquad a \ge0, $$ the inverse local time at position $x$.
 We define now  for bounded Borel functions $g$ with compact support,
\begin{equation}\label{def:Wtilde}
 W(g) := \int_{-\infty}^\infty g(L(t,X_t),X_t) {\rm sgn}(B_t) \d B_t.
\end{equation}

The stochastic integral  has to be understood as an integral with respect to the Brownian motion $B^{(r)}:=(B_{t+T_r},\, t\ge 0)$ where $r$ is any positive real such that $g(\ell,x)=0$ for all $x\ge r$ (that $B^{(r)}$ is a standard Brownian motion  comes from the fact that $(B'_{T_r+t},\, t\in [0,|T_r|])$ is distributed as $(B_{t},\, t\in [0,T_{-r}])$).  We will see in the following theorem that $W$ defines a white noise.  

Similarly to Notation \ref{n:<h}, for $x\in \r$, we can consider the process $(X^{-,x}_u,\,u\ge 0)$ obtained by gluing the excursions of $X$ below $x$:    that is we set for $t\in \r$, $A_t^{-,x} := \int_{-\infty}^t {\bf 1}_{\{X_s \le  x \}} {\rm d} s$, $\alpha_u^{-, x}  := \inf\{t\in \r:\, A_t^{-,x} > u\}$, and $X^{-,x}_u := X_{\alpha_u^{-,x}}$ for $u\ge 0$. The excursion filtration ${\mathcal E}^-$ is defined as  ${\mathcal E}^-_x:=\sigma(X_u^{-,x},\, u \ge0)$, for all $x\in \r$.

For ${\mathcal E}^+_x$,  we let $A_t^{+,x}:=t$ if $t\le T_x$, and for $t>T_x$, $A_t^{+,x} := T_x+ \int_{T_x}^t {\bf 1}_{\{X_s >  x \}} {\rm d} s$.  For $s\in \r$, we set $\alpha_s^{+, x}  := \inf\{t\in \r,\, A_t^{+,x} > s\}$, and  $ X^{+,x}_s := X_{\alpha_s^{+,x}}$.  The excursion filtration ${\mathcal E}^+$ is defined as   ${\mathcal E}^+_x :=  \sigma(X_s^{+,x},\, s \in \r)$,   for all $x\in \r$. Note that ${\mathcal E}^+_x$ is decreasing in $x$ whereas 
${\mathcal E}^-_x$ is increasing.

\begin{theorem}\label{t:mainbis}
Fix $\mu>0$. 
Equation \eqref{def:Wtilde} defines  a white noise on $\r_+\times \r$. 

(i)  Almost surely for all $a>0$, $r\in \r$ and $h\in [0,|I_{\tau_a^r}-r|]$,
\begin{equation}\label{eq:main1bis}
L(\tau_a^r,  r- h )
=
a - 2 \int_{r-h}^r  W\big([0, L(\tau_a^r, x)], \d x\big) + \left(2-{2\over \mu}\right) h, 
\end{equation}
\noindent where $ I_t:=\inf_{-\infty < s \le t} X_s,   t \in \r,$  denotes  the infimum process of $X$.

(ii)  Almost surely for all $a>0$, $r\in \r$ and $h\ge 0$, \begin{equation}\label{eq:main2bis}
L(\tau_a^r, r+h) 
=
 a+ 2 \int_r^{r+h}  W\big([0, L(\tau_a^r, x)], \d x\big)
+ {2\over \mu} h.
\end{equation}
\end{theorem}

The stochastic integral $\int_{r-h}^r  W\big([0, L(\tau_a^r, x)], \d x\big)$ is the stochastic integral with respect to the martingale measure $ M^{(r)}_h(\bullet):=W({\bf 1}_{\bullet\times [r-h,r]})$, $h\ge 0$,   associated to the filtration $({\mathcal E}_{r-h}^+)_{h\ge 0}$. The stochastic integral $\int_r^{r+h}  W\big([0, L(\tau_a^r, x)], \d x\big)$ is the stochastic integral with respect to the martingale measure $ \widehat M^{(r)}_h(\bullet):=W({\bf 1}_{\bullet\times [r,r+h]})$, $h\ge 0$, associated to the filtration $({\mathcal E}_{r+h}^-)_{h\ge 0}$.

\bigskip

\noindent {\it Proof}.
 Notice that for any $r\in \r$, $X^{(r)}:=(X_{T_r+t}-r,\,t\ge 0)$ is distributed as $(X_t,\,t\ge 0)$. Therefore we can apply Theorem \ref{t:main} to $X^{(r)}$.  As in \eqref{def:W}, we define $W^{(r)}$ the white noise associated to $X^{(r)}$: for any Borel function  $g:\r_+\times \r \to \r$ such that $\int_{\r_+} \d \ell \int_{\r} g^2(\ell,x) \d x<\infty$, $$W^{(r)}(g):= \int_0^\infty g\big(L^{(r)}(t, X^{(r)}_t), X^{(r)}_t\big)  {\rm sgn}(B^{(r)}_t) \d B^{(r)}_t,$$
 
\noindent where $L^{(r)}(\cdot, \cdot)$ denote the local times of $X^{(r)}$. Let $\tau^{x,(r)}_{a}$ be the associated inverse local times. By Theorem \ref{t:main} (i) applied to $X^{(r)}$,  for $a>0$,
\begin{equation}\label{eq:1bisproof}
L^{(r)}(\tau_a^{0,{(r)}},  - h )
=
a - 2 \int_{-h}^0 W^{(r)}\big([0, L^{(r)}(\tau_a^{0,(r)}, x)], \d x\big) + \left(2-{2\over \mu}\right) h, \qquad h\in [0,|I_{\tau_a^{0,(r)}}^{(r)}|],
\end{equation}
\noindent where $ I_t^{(r)}:=\inf_{0\le s \le t} X_s^{(r)},   t\ge0,$  is the infimum process of $X^{(r)}$. Notice that $\tau_a^{0,(r)}=\tau_a^r-T_r$,  $L^{(r)}(t,  x ) = L(t+T_r,r+x)$ for $x\le0$, and for $t\ge 0$, $I_t^{(r)} = I_{t+T_r} -r$ where  $I_t:=\inf_{-\infty < s \le t} X_s$. Moreover, for any $x\le y$,
$$
W^{(r)}({\bf 1}_{A\times [x,y]})
=
\int_0^\infty {\bf 1}_{A\times [x,y]}(L^{(r)}(t,X_t^{(r)}),X_t^{(r)}) {\rm sgn}(B_t^{(r)}) \d B_t^{(r)}
$$

\noindent where we recall that $B_t^{(r)}=B_{t+T_r}$, $t\ge 0$.  We deduce that if $x\le y \le 0$ and $A$ is a bounded Borel set of $\r_+$,
\begin{eqnarray}
W^{(r)}({\bf 1}_{A\times [x,y]})
&=& \nonumber
\int_{0}^\infty {\bf 1}_{A\times [x,y]}(L(t+T_r,X_{t+T_r}),X_{t+T_r}-r) \, {\rm sgn}(B_t^{(r)}) \d B_t^{(r)}
\\
&=& \nonumber
\int_{T_r}^\infty {\bf 1}_{A\times [r+x,r+y]}(L(t,X_{t}),X_{t})\, {\rm sgn}(B_t) \d B_t\\
&=&
W({\bf 1}_{A\times [r+x,r+y]}). \label{eq:wr}
\end{eqnarray} 

\noindent We deduce that \eqref{def:Wtilde} defines a white noise on $\r_+ \times (-\infty,r)$, hence on $\r_+\times \r$ since $r$ can be made arbitrary large. Equation \eqref{eq:1bisproof} becomes
$$
L(\tau_a^r,  r - h )
=
a - 2 \int_{r-h}^r W\big([0, L(\tau_a^r, x)], \d x\big) + \left(2-{2\over \mu}\right) h, \qquad h\in [0,|I_{\tau_a^r}-r|].
$$
 
\noindent It is \eqref{eq:main1bis}. We prove now \eqref{eq:main2bis}.  Let $r'>r+h$ arbitrary. Using Tanaka's formula  applied to $(X_{T_{r'}+t},\, t\ge 0)$, we have for any $x\in \r$ and $t\ge T_{r'}$,
\begin{equation}\label{eq:tanaka-r}
(X_t-x)^-=(r'-x)^- - \int_{T_{r'}}^t {\bf 1}_{\{ X_s \le x\}} \d X_s + \frac12 (L(t,x)-L(T_{r'},x)).
\end{equation}

\noindent Taking $t=\tau_a^r$,  and $x=r$ then $x=r+h$, we get
$$
L(\tau_a^r, r+h)
=
2h+ a+  2 \int_{T_{r'}}^{\tau_a^r} {\bf 1}_{\{r < X_s \le r+h\}} \d X_s
=
a+  2 \int_{T_{r'}}^{\tau_a^r} {\bf 1}_{\{r < X_s \le r+h\}} {\rm sgn}(B_s)\d B_s +{2\over \mu} h
$$

\noindent where the second equality follows from \eqref{dXs}   and the fact that $\int_{T_{r'}}^{\tau_a^r} {\bf 1}_{\{r < X_s \le r+h\}}   \d I_s =\int_{T_{r'}}^{T_r} {\bf 1}_{\{I_s \le r+h\}}   \d I_s = -h$.  By Proposition \ref{p:main} applied to $(X^{(r')},B^{(r')})$, 
\begin{eqnarray*}
\int_{T_{r'}}^{\tau_a^r} {\bf 1}_{\{r < X_s \le r+h\}} {\rm sgn}(B_s)\d B_s
&=&
\int_0^{\tau_a^{r-r',(r')}} {\bf 1}_{\{r-r' < X_s^{(r')} \le r+h-r'\}} {\rm sgn}(B_s^{(r')})\d B_s^{(r')}
\\
&=&
\int_{r-r'}^{r+h-r'}   W^{(r')}\big([0, L(\tau_a^r, x)], \d x\big)\\
&=&
\int_{r}^{r+h}   W\big([0, L(\tau_a^r, x)], \d x\big)
\end{eqnarray*}

\noindent by \eqref{eq:wr} applied to $r=r'$. This proves \eqref{eq:main2bis} and completes the proof of the theorem. $\Box$

\end{document}